\documentclass[manuscript]{copernicus}

\usepackage{amsmath,amssymb,amsfonts,amsthm,mathtools}
\usepackage[dvipsnames]{xcolor}
\usepackage[toc,page]{appendix}
\usepackage{csl}

\usepackage{subcaption}

\def\!#1{\mathcal{#1}}
\def\*#1{\boldsymbol{\mathbf{#1}}}
\def\|#1{\textnormal{#1}}
\def\##1{\mathfrak{#1}}

\def\norm#1{\left\lVert#1\right\rVert}

\DeclareMathOperator*{\argmin}{arg\,min}
\DeclareMathOperator*{\argmax}{arg\,max}
\DeclareMathOperator{\tr}{tr}

\newcommand{\xfrv}{X^\|f}
\newcommand{\xarv}{X^\|a}
\newcommand{\xhfrv}{\hat{X}^\|f}
\newcommand{\xharv}{\hat{X}^\|a}
\newcommand{\xf}{\*X^\|f}
\newcommand{\xa}{\*X^\|a}

\newcommand{\ensNf}{N^{\|f}}
\newcommand{\ensNa}{N^{\|a}}

\graphicspath{ {Figures/} }

\makeatletter
\renewcommand*\env@matrix[1][c]{\hskip -\arraycolsep
    \let\@ifnextchar\new@ifnextchar
    \array{*\c@MaxMatrixCols #1}}
\makeatother

\usepackage{amsmath}
\usepackage{paralist}
\usepackage{graphics} %% add this and next lines if pictures should be in esp format
\usepackage{epsfig} %For pictures: screened artwork should be set up with an 85 or 100 line screen
\usepackage{graphicx}  \usepackage{epstopdf}%This is to transfer .eps figure to .pdf figure; please compile your paper using PDFLeTex or PDFTeXify.
\usepackage{cleveref}
\theoremstyle{definition}

\newtheorem{remark}{Remark}

\csltitle{A Stochastic Covariance Shrinkage Approach to Particle Rejuvenation in the Ensemble Transform Particle Filter}
\cslauthor{Andrey A Popov, Amit N Subrahmanya, Adrian Sandu}
\cslyear{20}
\cslreportnumber{4}
\cslemail{apopov@vt.edu, amitns@vt.edu, sandu@cs.vt.edu}
\begin{document}
    \nolinenumbers
    
    \title{A Stochastic Covariance Shrinkage Approach to Particle Rejuvenation in the Ensemble Transform Particle Filter}
    
    \Author[1]{Andrey A}{Popov}
    \Author[1]{Amit N}{Subrahmanya}
    \Author[1]{Adrian}{Sandu}
    
    \affil[1]{Computational Science Laboratory, Department of Computer Science, Virginia Tech, 2202 Kraft Drive, Blacksburg, VA, 24060, USA}

    %% The [] brackets identify the author with the corresponding affiliation. 1, 2, 3, etc. should be inserted.
    
    %% If an author is deceased, please mark the respective author name(s) with a dagger, e.g. "\Author[2,$\dag$]{Anton}{Smith}", and add a further "\affil[$\dag$]{deceased, 1 July 2019}".
    
    %% If authors contributed equally, please mark the respective author names with an asterisk, e.g. "\Author[2,*]{Anton}{Smith}" and "\Author[3,*]{Bradley}{Miller}" and add a further affiliation: "\affil[*]{These authors contributed equally to this work.}".

    \correspondence{Andrey A Popov (apopov@vt.edu)}
    
    \runningtitle{Foresight ETPF}
    
    \runningauthor{AA Popov, AN Subramanya, A Sandu}

    \received{}
    \pubdiscuss{} %% only important for two-stage journals
    \revised{}
    \accepted{}
    \published{}

    \csltitlepage
    % Enter the first author's name and address:
    %\centerline{\scshape Andrey A. Popov$^*$, Amit N. Subrahmanya and Adrian Sandu}
    %\medskip
    %{\footnotesize
    %    % please put the address of the first author
    %    \centerline{ 2202 Kraft Drive}
    %    \centerline{ Blacksburg, VA 24060, USA}
    %} % Do not forget to end the {\footnotesize by the sign }
    
    \bigskip
    
    % The name of the associate editor will be entered by an editorial staff
    % "Communicated by the associate editor name" is not needed for special issue.
    %\centerline{(Communicated by the associate editor name)}

    \maketitle
    
    \begin{abstract}
        Rejuvenation in particle filters is necessary to prevent the collapse of the weights when the number of particles is insufficient to sample the high probability regions of the state space. Rejuvenation is often implemented in a heuristic manner by the addition of stochastic samples that widen the support of the ensemble. This work aims at improving canonical rejuvenation methodology by the introduction of additional prior information obtained from climatological samples; the dynamical particles used for importance sampling are augmented with samples obtained from stochastic covariance shrinkage. The ensemble transport particle filter, and its second order variant, are extended with the proposed rejuvenation approach. Numerical experiments show that modified filters significantly improve the analyses for low dynamical ensemble sizes.
    \end{abstract}

    %%%%%%%%%%%%%%%%%%%%%
    \section{Introduction}
    %%%%%%%%%%%%%%%%%%%%%
    
    Ensemble-based data assimilation~\citep{asch2016data,law2015data,reich2015probabilistic} aims to estimate our uncertainty about the state of some dynamical system through an ensemble of possible states. The assumed underlying distribution of these states is taken to be our uncertainty.
    
    Oftentimes the ensemble cannot describe the distribution of our uncertainty to sufficient accuracy. The curse of dimensionality~\citep{karpatne2018book} ensures that true descriptions of general high-dimensional distributions stay out of our reach. Several techniques such as the principal of maximum entropy~\citep{jaynes2003probability} exist in order to attempt to alleviate some of the burden by prescribing a distribution to a set of data by constraining the choice with  known or estimated quantities and qualities. For instance, the ensemble Kalman filter~\citep{burgers1998analysis,evensen1994sequential,evensen2009data}, could be thought of as an abuse of the principle of maximum entropy: by discarding information about the underlying dynamical system, and assuming that the ensemble only gives information about a mean (that lives in $\mathbb{R}^n$), and a covariance, the underlying distribution of our uncertainty is taken to be normal. In such a way, any application of Bayes' rule has to transform our assumed prior normal distribution into our assumed posterior normal distribution. 
    
    Previous work~\citep{popov2020stochastic} has focused on augmenting the information represented by the ensemble with information derived from covariance shrinkage through a surrogate ensemble in the ensemble transport Kalman filter. In this paper, we extend this idea to the ensemble transport particle filter (ETPF)~\citep{reich2013nonparametric}. The ETPF attempts to use importance sampling~\citep{liu2008monte}, not for a solution to the problem of Bayesian inference, but to simply transport a given ensemble into one that is equally sampled from some distribution whose moments, in the ensemble limit, approach the moments of the correct posterior distribution. This means that the underlying distribution from which our posterior is sampled, could be potentially very far, in information distance, to the optimal posterior distribution, for a finite ensemble.
    
   This work explores a new approach to particle rejuvenation, which is necessary to prevent weight collapse in particle filters. Instead of heuristics, the approach makes use of prior information to enrich the ensemble subspace. Our contributions are as follows: we introduce an alternative way of performing particle rejuvenation in ETPF by incorporating climatological covariance information. We accomplish this by augmenting the dynamical (model) ensemble with synthetic anomalies with optimal scaling,  accompanied by a statistically correct estimator. We show that this method of performing particle rejuvenation significantly improves the analysis RMSE for low dynamical ensemble sizes.
    
    This paper is organized as follows.  \Cref{sec:bayesianinference} reviews the concept of Bayesian inference with the addition of prior information, and its use in importance sampling. \Cref{sec:ETPF} introduces the ensemble transform particle filter and its canonical rejuvenation heuristic. The concept of stochastic covariance shrinkage is proposed in \Cref{sec:stochasticshrinkage}, and ETPF is extended to make use of this shrinkage. Numerical experiment results are reported in \Cref{sec:numerical}. Finally, concluding remarks are drawn in \Cref{sec:conclusions}.

    %%%%%%%%%%%%%%%%%%%%%
    \section{Optimal coupling with prior information and the ensemble transform particle filter}
    \label{sec:bayesianinference}
    %%%%%%%%%%%%%%%%%%%%%
    
    Bayesian inference~\citep{jaynes2003probability} aims at transforming prior information about the state of a system (represented by the distribution of a random variable $\xfrv$), additional qualitative and quantitative information ($P$), which we will use to stand for , 
    and information obtained by observing the system ($Y$), into combined posterior information ($\xarv$):
    \begin{equation}
        \label{eqn:Bayes}
        \pi_{\xarv}(X \mid Y, P) = \frac{\pi_{Y | X} (Y \mid X, P)\,\pi_{\xfrv}(X \mid P)}{\pi_Y(Y \mid P)},
    \end{equation}
    where $\pi_{\xfrv}(x\mid p)$ represents the prior state probability density conditioned by all other relevant information, and $\pi_{Y | \xfrv} (y\mid x, p)$ is the observational likelihood conditioned by the the forecast and the prior information.
    Here we consider the finite dimensional case where $\xfrv,\xarv \in \mathbb{R}^n$, $Y \in \mathbb{R}^{m}$, and the supports of the  probability densities $\pi_{\xfrv}$ and $\pi_{\xarv}$ are subsets of the respective spaces.
    
    Classical particle filtering~\citep{reich2015probabilistic} represents state distributions by collections of particles, i.e., ensembles of samples. 
    Specifically, consider an ensemble of $\ensNf$ particles $\xf = [\xf_1,\dots \xf_{\ensNf}] \in \mathbb{R}^{n \times \ensNf}$. The prior distribution density $\pi_{\xfrv}$ is approximated weakly by the corresponding empirical measure
    \begin{equation} 
        \label{eqn:prior-measure-distribution}
        \pi_{\xhfrv}(X\mid P) = \sum_{i=j}^{\ensNf} w^\|f_j\, \delta_{X-\xf_j},
    \end{equation}
    where $w^\|f_j$ for $1\leq j\leq \ensNf$ are the prior importance weights associated with each particle. Similarly, the posterior density is approximated weakly by an empirical measure based on the same sample values (particle states) but with different posterior importance weights  $w^\|a_j$ for $1\leq j\leq \ensNf$:
    \begin{equation} 
        \label{eqn:posteriot-measure-distribution}
        \pi_{\xharv}(x\mid y,p) = \sum_{i=j}^{\ensNf} w^\|a_j\, \delta_{x-\xf_j},
    \end{equation}
    The posterior importance sampling weights are obtained from \cref{eqn:Bayes}:
    \begin{equation}
        \label{eqn:posterior-weights-IS}
        w^\|a_j \propto  \pi_Y(y \mid \xfrv_j, P)\,\pi_P(P \mid \xfrv_j) = w^\|f_j\, \pi_Y(y \mid \xfrv_j, P).
    \end{equation}

    The ensemble of weights is denoted by $\*w = [w_1,\dots w_{\ensNf}]^T$, and $\*w^\|f$ and $\*w^\|a$ refer to the forecast and analysis weights respectively. Using \eqref{eqn:posteriot-measure-distribution} and \eqref{eqn:posterior-weights-IS} empirical estimates of the posterior mean and covariance,
    \begin{equation}\label{eqn:mean-cov-estimate}
        \begin{gathered}
            \bar{\*x}^\|a = \sum_{j=1}^{\ensNf} w^\|a_j \xf_j,\\
            \*\Sigma_{X^\|a} = \frac{\ensNf}{\ensNf-1}\xf_j \left(\text{diag}\left(\*w^\|a\right) - \*w^\|a \*w^{\|a,T}\right){\*X^{\|f,T}_j},
        \end{gathered}
    \end{equation}
    respectively, are obtained by the importance sampling approach~\citep{liu2008monte}.

    The goal of particle filtering (with resampling) is to find an ensemble $\xa\in \mathbb{R}^{n\times\ensNa}$ of $\ensNa$ realizations of the random variable $\xarv$ that represents the posterior distribution $\pi_{\xarv}$ with equal weights. 
    Specifically, the the posterior density is approximated weakly by the empirical measure 
    \begin{equation} 
        \label{eqn:posterior-pdf-IS}
        \pi_{\xharv}(x\mid p) = \sum_{j=1}^{\ensNa} \frac{1}{\ensNa} \delta_{x-\xa_j},
    \end{equation}
    where the importance sampling weights are uniform and equal to $1/\ensNa$ (so as to be equally likely). 
    We impose that the empirical mean \eqref{eqn:mean-cov-estimate} is preserved by \eqref{eqn:posterior-pdf-IS}:
    \begin{equation}
        \bar{\*x}^\|a = \sum_{j=1}^{\ensNa} \frac{1}{\ensNa} \xa_j = \sum_{j=1}^N w^\|a_j \xf_j =  \xf \*w^\|a.
    \end{equation}

    The optimal coupling~\citep{mccann1995existence,reich2015probabilistic} between the prior empirical distribution \cref{eqn:prior-measure-distribution} and the posterior empirical distribution \cref{eqn:posterior-pdf-IS}, can be defined as an ensemble transformation,
    \begin{equation}\label{eqn:optimal-transport}
        \xa = \xf\,\*T^*,
    \end{equation}
    where $\*T^* \in \mathbb{R}^{\ensNf\times\ensNa}$ is the solution to the optimal transport Monge-Kantorovich problem~\citep{villani2003topics}. 
    The discrete optimal transportation problem is
    \begin{equation}\label{eqn:mongeproblem}
        \begin{gathered}
            \*T^* = \argmin_{\*T} \sum_{\substack{1\leq j\leq \ensNf\\ 1\leq k\leq \ensNa}} \*T_{j,k}\norm{\xf_j - \xf_k}_2^2\\ 
            \text{subject to}\quad \*T\*1_{\ensNa} = \ensNa\*w^\|a,\, \*T^T\*1_{\ensNf} = \*1_{\ensNf},\, \*T_{i,j} \geq 0,
        \end{gathered}
    \end{equation}
    where the distance measure of squared Euclidean distance is taken for a provably unique solution to the Monge-Kantorovich problem to exist~\citep{mccann2011five}. The vector of ones of size $q$ is represented by $\*1_q$.
    The problem~\cref{eqn:mongeproblem} is a linear programming problem.
    
    The discrete optimal transport \cref{eqn:optimal-transport} formulation begets a mapping $\xa = \Psi_{\ensNf, \ensNa}(\xf)$, that, in the limit of ensemble size ($\ensNf = \ensNa \to \infty$), converges weakly to a mapping $\Psi$, such that $\xarv = \Psi(\xfrv)$, which has the exact desired distribution given by \cref{eqn:Bayes} \cite[Theorem 5.19 ]{reich2015probabilistic}. We believe, but don't prove, that this is likely when $N^\|f$ and $N^\|a$, are not equal.

    The standard ETPF~\citep{reich2013nonparametric} makes the assumption that the prior and posterior ensemble sizes are the same, $N \coloneqq \ensNa = \ensNf$ in \eqref{eqn:mongeproblem}. A second order extension to the ETPF (which we will call ETPF2 here)~\citep{acevedo2017second} modifies the optimal transport equation \eqref{eqn:optimal-transport} as follows: 
    \begin{equation}
        \xa = \xf\left(\*T^* + \*D\right),
    \end{equation}
    where the additional term $\*D$ is a matrix that ensures that the empirical covariance estimate $\*\Sigma_{X^\|a}$ from \eqref{eqn:mean-cov-estimate} is preserved by \eqref{eqn:posterior-pdf-IS}.

    %%%%%%%%%%%%%%%%%%%%%
    \section{Particle Rejuvenation in ETPF}
    \label{sec:ETPF}
    %%%%%%%%%%%%%%%%%%%%%
    
    {Particle and ensemble-based filters often underrepresent uncertainty~\citep{asch2016data} due to the relatively small number of samples when compared to the dimension of the state and data spaces. Over several data assimilation cycles multiple particle start carrying either unimportant or redundant information, which leads to weight collapse or to ensemble degeneracy ~\citep{strogatz2018nonlinear}. To alleviate these effects, methods such as inflation~\citep{anderson2001ensemble,popov2020explicit}, rejuvenation~\citep{reich2013nonparametric}, and resampling~\citep{reich2015probabilistic} have been developed.
    
    In order to avoid ensemble collapse, ETPF employs a particle rejuvenation approach ~\citep{acevedo2017second,reich2013nonparametric,chustagulprom2016hybrid} that perturbs the analysis ensemble by a random sampling from the ensemble of prior anomalies, 
    \begin{equation}\label{eqn:naive-rejuvenation}
        \xa \leftarrow \xa + \sqrt{\frac{\tau}{N - 1}}\*A^\|f\*\eta\left(\*I_N - N^{-1}\*1_N\*1_N^T\right),
    \end{equation}
    where $\*\eta \sim \left( \mathcal{N}(0,1) \right)_{N \times N}$ is a matrix of i.i.d.  samples from the standard normal distribution of size $N$, the factor $\tau$ is treated as a hyperparameter that controls the magnitude of the correction, and the ensemble anomalies
    \begin{equation}\label{eqn:anomalies}
        \*A^\|f = \xf\left(\*I_N - N^{-1}\*1_N\*1_N^T\right),
    \end{equation}
    are defined as the ensemble of deviations from the sample mean. Of note is the fact that the extra term $\left(\*I_N - N^{-1}\*1_N\*1_N^T\right)$ in \eqref{eqn:naive-rejuvenation} ensures that the introduction of the random matrix $\*\eta$ does not modify the mean of $\*X^\|a$. 
    This is due to the fact that,
    \begin{equation}\label{eq:consistency}
        \begin{gathered}
            \left(\*I_N - N^{-1}\*1_N\*1_N^T\right)\*1_N = \*0_N,\\
            \*1_N^T\left(\*I_N - N^{-1}\*1_N\*1_N^T\right) = \*0_N^T.
        \end{gathered}
    \end{equation}

    Notice that if we define the matrix,
    \begin{equation}\label{eqn:canonicalrejuv}
        \*B \coloneqq \sqrt{\frac{\tau}{N - 1}}\left(\*I_N - N^{-1}\*1_N\*1_N^T\right)\*\eta\left(\*I_N - N^{-1}\*1_N\*1_N^T\right),
    \end{equation}
    it is possible to write the ETPF with rejuvenation \eqref{eqn:naive-rejuvenation} as,
    \begin{equation}\label{eqn:perturb-ensemble}
        \begin{aligned}
            \*X^{\|a} &= \*X^{\|f}\*T^* + \sqrt{\frac{\tau}{N-1}}\*A^{\|f}\*\eta\left(\*I_N - N^{-1}\*1_N\*1_N^T\right)\\
            &= \*X^{\|f}\*T^* + \sqrt{\frac{\tau}{N-1}}\*X^{\|f}\left(\*I_N - N^{-1}\*1_N\*1_N^T\right)\*\eta\left(\*I_N - N^{-1}\*1_N\*1_N^T\right)               \\
            &=\*X^{\|f}\left[\*T^* + \sqrt{\frac{\tau}{N-1}}\left(\*I_N - N^{-1}\*1_N\*1_N^T\right)\*\eta\left(\*I_N - N^{-1}\*1_N\*1_N^T\right)\right]               \\
            &= \*X^{\|f}\, \widetilde{\*T} \quad \textnormal{with} \quad \widetilde{\*T} \coloneqq \*T^\ast + \*B,
        \end{aligned}
    \end{equation}
    with the matrix $\*B$ acting as a stochastic perturbation of the optimal transport operator $\*T^*$, such that $\*B$ preserves the constraints $\*1_N^T \*B = \*0_N^T$ and $\*B \*1_N = 
    \*0_N$ in \cref{eqn:mongeproblem} because of the results in \eqref{eq:consistency}. This, of course, immediately calls into question the optimality of the transport for a finite ensemble, as adding this type of $\*B$ matrix is perturbing the transport mapping $\widetilde{\*T}$ away from the optimum $\*T^\ast$.

    %%%%%%%%%%%%%%%%%%%%%
    \section{Particle Rejuvenation Through Stochastic Shrinkage}
    \label{sec:stochasticshrinkage}
    %%%%%%%%%%%%%%%%%%%%%
    
    In the context of ensemble methods, covariance shrinkage~\citep{Sandu_2019_Covariance-parallel,Sandu_2015_covarianceShrinkage,Sandu_2014_EnKF_SMF} is used, similar to other canonical covariance tapering techniques such as  inflation~\citep{anderson2001ensemble, popov2020explicit}, localization~\citep{anderson2012localization,hunt2007efficient,Sandu_2017_Covariance-Cholesky, Sandu_2015_SCALA,petrie2008localization,zhang2010ensemble}, to enrich the information represented by an undersampled covariance matrix.
    
%An ensemble estimator of the mean would have $n$ degrees of freedom, while a covariance estimator would have $\!O(n^2)$ degrees of freedom. 
%        If an estimate of the covariance already exists, then the problem of determining the covariance is reduced to determining the deviation from said estimate. 
%        Covariance shrinkage is a statistical attempt to formalize this concept. 
        From a Bayesian perspective, covariance shrinkage seeks to incorporate additional prior information on error correlations into the analysis, in order to enhance the inference.
        In many data assimilation models, climatological covariance information is often available, i.e., it is known prior information.
        Climatological covariances are typically precomputed or derived from climatological models and are often employed in variational data assimilation \cite{lorenc2015comparison}.
    
    Following~\citep{popov2020stochastic}, we describe the stochastic covariance shrinkage technique.
    Instead of perturbing the transform matrix as in \eqref{eqn:perturb-ensemble}, we instead consider enhancing the dynamic ensemble $\*X^\|f \in \mathbb{R}^{n\times N}$ with an $M$ member synthetic ensemble $\!X^\|f  \in \mathbb{R}^{n\times M}$ of samples independent of the dynamical ensemble. 
    This leads to the total $N+M$ members ensemble:
    \begin{equation}\label{eqn:total-ensemble}
        \#X^\|f = \begin{bmatrix}\*X^\|f & \!X^\|f\end{bmatrix}  \in \mathbb{R}^{n \times (N + M)}
    \end{equation}
    with weights $\#w^\|f = [w^\|f_1, \dots w^\|f_{N+M}]$. We make the ansatz that each ensemble member is distributed as 
    \begin{equation}\label{eq:full-dist-assumption}
        \#X^\|f_{:,i} \sim \pi_{X^\|f}(X \mid P),
    \end{equation}
    which is the full distribution of the forecast conditioned by the prior information that we have provided to the algorithm. Note that \eqref{eq:full-dist-assumption} is not the empirical measure distribution \eqref{eqn:prior-measure-distribution},
    that only has information from the ensemble members, but rather the `exact' distribution that is assumed to contain all the information from the forecast. THis is because we are now incorporating more prior information $P$ in the form of climatological information.
    
    Taking \eqref{eqn:anomalies} to be the anomalies of the dynamic ensemble, and 
    \begin{equation}\label{eqn:anomalies-static}
        \!A^\|f = \!X^\|f\left(\*I_M - M^{-1}\*1_M\*1_M^T\right),
    \end{equation}
    to be the anomalies of the synthetic ensemble, the total empirical (unbiased) covariance can be written as,
    \begin{equation}
        \*\Sigma_{\#X^\|f} = \*\Sigma_{\*X^\|f} + \*\Sigma_{\!X^\|f},
    \end{equation}
    where the constituent covariances are defined in terms of the weights
    \begin{equation}
        \*\Sigma_{\*X^\|f} = \sum_{i=1}^N w^\|f_i\frac{N}{N-1} \*A_{:,i}^\|f\*A_{:,i}^{\|f,T},\quad \*\Sigma_{\!X^\|f} = \sum_{i=1}^M w^\|f_{N+i}\frac{M}{M-1}  \!A_{:,i}^\|f\!A_{:,i}^{\|f,T}.
    \end{equation}

    In the covariance shrinkage approach, the sample mean of the synthetic ensemble is assumed to be the sample mean of the dynamic ensemble,
    \begin{equation}
        \!X^\|f\*1_M = \*X^\|f\*1_N,
    \end{equation}
    by construction and \eqref{eq:consistency}, thus requiring that only the synthetic ensemble anomalies need to be determined. Taking a covariance matrix $\!P$ that represents a climatological estimate of the covariance, we sample the anomalies from some unbiased distribution with a scaling, by a factor $\mu$, of said covariance. In the Gaussian case
    \begin{equation}
        \!A^\|f_{:,i} \sim \!N(\*0_n, \mu\!P), 
    \end{equation}
    where $\mu$ is a scaling factor defined later. An alternate choice of distribution that we explore is the symmetric Laplace distribution~\citep{kozubowski2013multivariate},
    \begin{equation}
        \!A^\|f_{:,i} \sim \!L(\*0_n, \mu\!P), 
    \end{equation}
    which is described by the pdf
    \begin{equation}
        \pi(\*x) = \frac{2}{\left[(2\pi)^n \mu\!P\right]^{\frac{1}{2}}} \left(\frac{\*x^T \!P^{-1} \*x}{2 \mu}\right)^{\frac{2 - n}{4}} K_{\frac{2 - n}{2}}\left(\sqrt{\frac{2\*x^T \!P^{-1} \*x}{\mu}}\right),
    \end{equation}
    where $K$ is the modified Bessel function of the second kind~\citep{olver2010nist}. The choice of Laplace distribution is motivated by robust statistics techniques~\citep{rao2017robust}. The resulting sampled covariance would therefore be an estimate of the scaled climatological covariance,
    \begin{equation}
        \*\Sigma_{\!X^\|f} \approx \mu\!P.
    \end{equation}

    \begin{remark}
        In order to stay consistent with the mean estimate, the anomalies are replaced with their sample mean zero counterparts
        \begin{equation}
            \!A^\|f \leftarrow  \!A^\|f\left(\*I_M - M^{-1}\*1_M\*1_M^T\right).
        \end{equation}
    \end{remark}
    
    The weights $\#w^\|f$ are divided into two classes: those that are associated with the dynamic ensemble, and those that are associated with the synthetic ensemble,
    \begin{equation}
        w^\|f_i = \begin{cases}
            1-\gamma & 1 \leq i \leq N\\
            \gamma & N+1 \leq i \leq N+M
        \end{cases},
    \end{equation}
    where the parameter $\gamma$ is known as the covariance shrinkage factor.
    
    One choice to calculate $\gamma$ is the Rao-Blackwell Ledoit-Wolf (RBLW) estimator~\citep{chen2009shrinkage} \cite[equation~(9)]{Sandu_2017_Covariance-Cholesky}, 
    \begin{equation}\label{eq:shr-rblw}
        \gamma_{\text{RBLW}} = \left[\vcenter{\hbox{$\displaystyle\frac{N - 2}{N(N+2)} + \frac{(n + 1)N - 2}{\hat{U}(\!P,\*\Sigma_{\*X^\|f})N(N+2)(n-1)}$}}\right] \land 1,
    \end{equation}
    where the sphericity factor,
    \begin{equation}\label{eqn:sphericity}
        \hat{U}(\!P,\*\Sigma_{\*X^\|f}) := \frac{1}{n - 1}\left(\frac{n \tr(\!C_i^2)}{\tr^2(\!C_i)} - 1\right),
        \quad\textnormal{with}\quad \!C \coloneqq \!P^{-\frac{1}{2}}\*\Sigma_{\*X^\|f} \!P^{-\frac{1}{2}},
    \end{equation}
    represents the mismatch between the climatological covariance (called the ``target'' in statistical literature) and sample covariance matrices.
    Note that if  $\!P = \*\Sigma_{\*X^\|f}$, then $\hat{U}(\!P,\*\Sigma_{\*X^\|f}) = 0$.
    In such a framework the scaling parameter for the climatological covariance is defined to be
    \begin{equation}
        \mu = \frac{\tr(\!C)}{n}\label{eqn:mu}.
    \end{equation}

    \begin{remark}
        The RBLW estimator \eqref{eq:shr-rblw} makes the assumption that the underlying distribution of the dynamic ensemble is Gaussian. Typically this assumption is violated for dynamical systems of interest.
    \end{remark}
    
    \begin{remark}
        In statistical literature, the target covariance is often taken to be the identity, $\!P = \*I$, which implies that $\!C = \*\Sigma_{\*X^\|f}$ in \eqref{eqn:sphericity}. The assumption that the target is a climatological covariance is natural generalization in the specific context of sequential data assimilation.
    \end{remark}
    
    \begin{remark}
        The scaling of the target matrix $\!P$ is not of any consequence. Let $\widetilde{\!P} = \beta\!P$ be a scalar scaling of the target matrix, then $\widetilde{\!C}_i = \frac{1}{\beta}\!C_i$, implying that $\widetilde{\mu}_i = \frac{1}{\beta}\mu_i$, rendering the matrix scaling inconsequential.
    \end{remark}

    The resulting analysis ensemble based on prior states and importance sampling weights of $N+M$ states is transported into an equally weighted posterior ensemble of $N$ states through the transformation
    \begin{equation}
        \xa = \#X^\|f \, \*T^*,
    \end{equation}
    where the optimal transport matrix $\*T^* \in \*R^{(N + M) \times N}$ is computed by solving \eqref{eqn:mongeproblem}.

    Recall that in the traditional method of rejuvenation \eqref{eqn:perturb-ensemble}, the optimal transport matrix is perturbed randomly into a nearby transport matrix; no new prior information is introduced. 
        We take a fundamentally different approach by incorporating ``unseen'' prior information derived from a climatological covariance.
        To this end, we augment of the empirical measure distribution \eqref{eqn:prior-measure-distribution} with samples from the climatological distribution, to accommodate the total ensemble \eqref{eqn:total-ensemble},
        \begin{equation} 
            \label{eqn:augmented-measure-distribution}
            \pi_{\xhfrv}(x\mid p) = \sum_{i=j}^{\ensNf} w^\|f_j\, \delta_{x-\#X^{\|f}_j},
        \end{equation}
        before the Monge-Kantorovich problem \eqref{eqn:mongeproblem} is solved. 
        In effect we are able to avoid ensemble collapse by enhancing the empirical measure distribution \eqref{eqn:augmented-measure-distribution} with new prior information, as opposed to a reweighing of the old prior information.
        We will denote our method as FETPF, standing for `foresight' ETPF, as we believe including prior information in the analysis procedure is a type of foresight.
    
    %%%%%%%%%%%%%%%%%%%%%%
    \subsection{Multiple Climatological Covariance Matrices}
    \label{sec:manyP}
    %%%%%%%%%%%%%%%%%%%%%%
    
    It is conceivable that multiple climatological models give rise to multiple climatological covariances, or alternatively multiple candidates for the most `common' behavior of the model is to be chosen.
    
    Given a collection of target covariances, $\{\!P_j\}_{j\in\!J}$, we must choose the appropriate covariance from which to sample. We consider the sphericity of the mismatch between the target and forecast covariances \cref{eqn:sphericity}.  Based on authors' numerical experience, we select the target covariance that corresponds to the highest sphericity of the mismatch:
    \begin{equation}
        \!P^* = \argmax_{\!P_j} \; \hat{U}(\!P_j,\*\Sigma_{\*X^\|f}),
    \end{equation}
    We can justify this choice by realizing that the smaller the sphericity, the closer our samples are to that of canonical rejuvenation techniques. The aim of climatological shrinkage is to introduce unknown information into our inference procedure, thus the target covariance with the highest mismatch introduces the highest amount of outside information.
    
    \begin{remark}
        It is also possible to construct `multi-target' shrinkage estimators~\citep{lancewicki2014multi} that consider all target matrices simultaneously.
    \end{remark}
    
    %%%%%%%%%%%%%%%%%%%%%
    \section{Numerical Experiments}
    \label{sec:numerical}
    %%%%%%%%%%%%%%%%%%%%%
    
    \begin{figure}[t]
        \centering
        \includegraphics[width=0.5\linewidth]{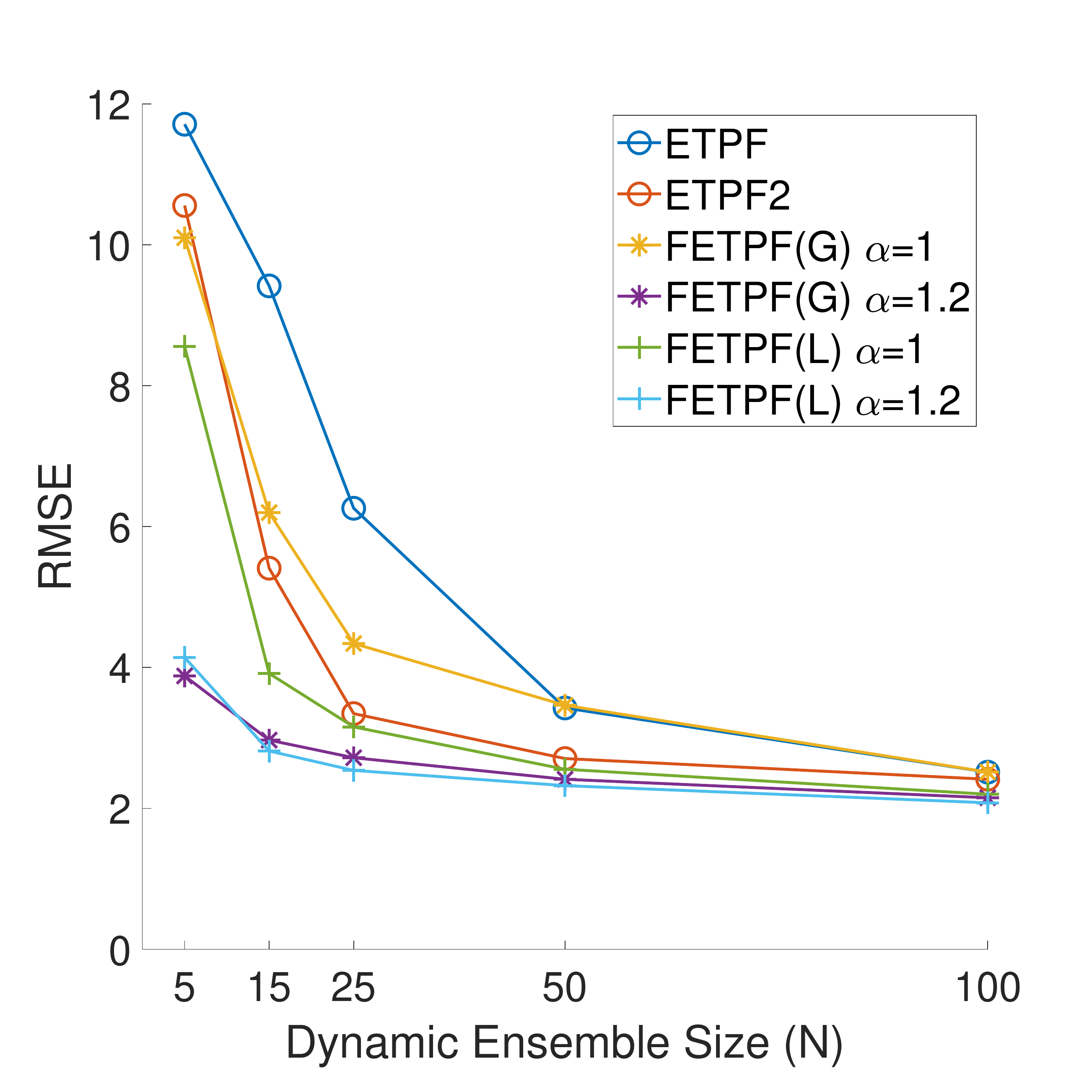} 
        \caption{Analysis RMSE versus dynamic ensemble size ($N$) of the Gaussian (G) and Laplace (L) covariance shrinkage approaches (FETPF) to particle rejuvenation for a synthetic ensemble size of $M=100$ and for various values of synthetic anomaly inflation, $\alpha$, from \eqref{eqn:synthetic-ensemble-inflation}, with respect to a canonical particle rejuvenation approach for first order ETPF and second order ETPF (denoted ETPF2)  for an optimal rejuvenation factor $\tau=0.04$. The target covariance is taken to be \eqref{eqn:cov-P}.}
        \label{fig:exp1}
    \end{figure}
    
    \begin{figure}[t]
        \centering
        \includegraphics[width=0.5\linewidth]{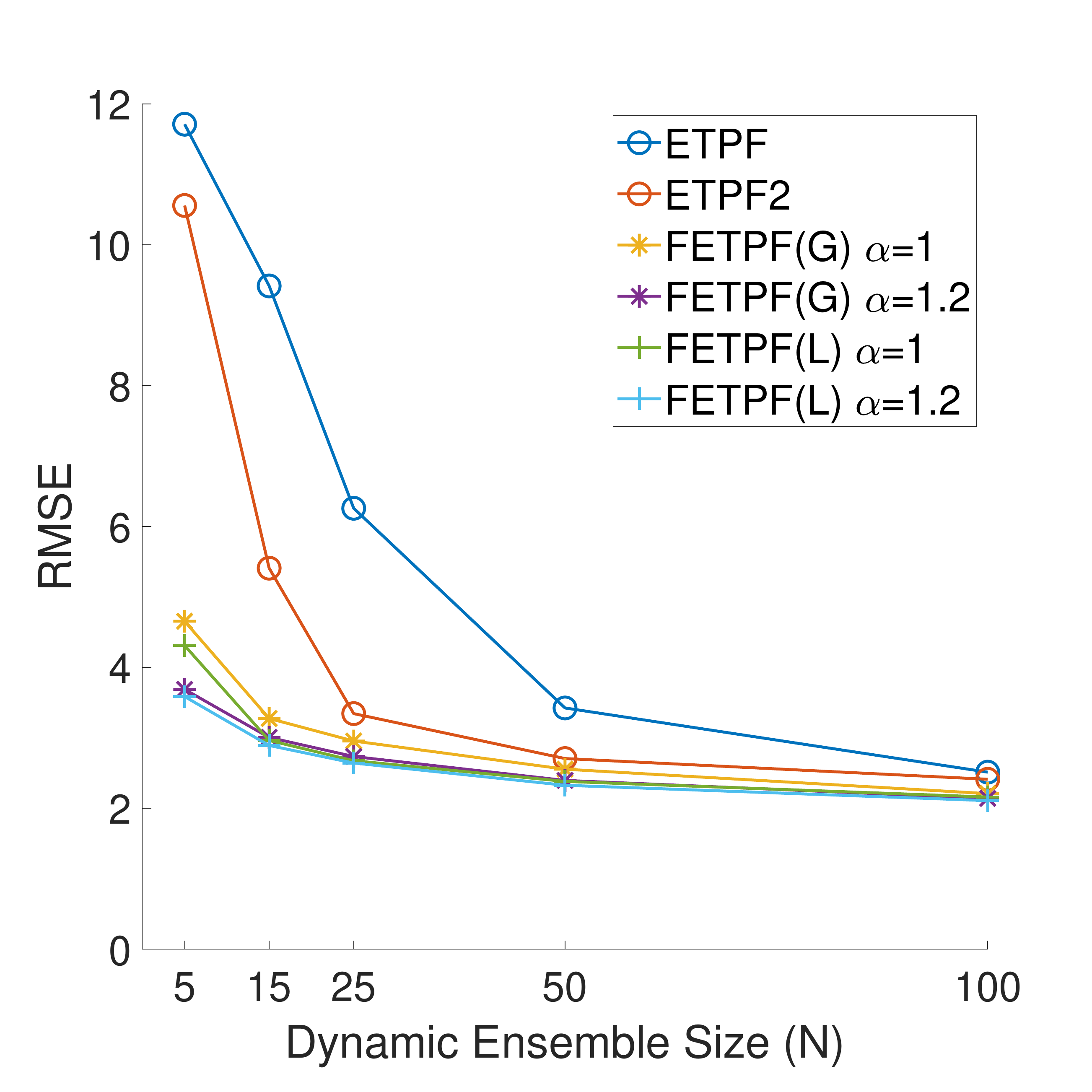} 
        \caption{Analysis RMSE versus Dynamic ensemble size ($N$) of the Gaussian (G) and Laplace (L) covariance shrinkage approaches (FETPF) to particle rejuvenation with the multi-target covariances \eqref{eqn:cov-P1-P2} for a synthetic ensemble size of $M=100$ and for various values of synthetic anomaly inflation ($\alpha$) with respect to a canonical particle rejuvenation approach for first order ETPF and second order ETPF (denoted ETPF2) for an optimal rejuvenation factor $\tau=0.04$.}
        \label{fig:exp2}
    \end{figure}
    
    \begin{figure}[t]
        \centering
        \includegraphics[width=0.5\linewidth]{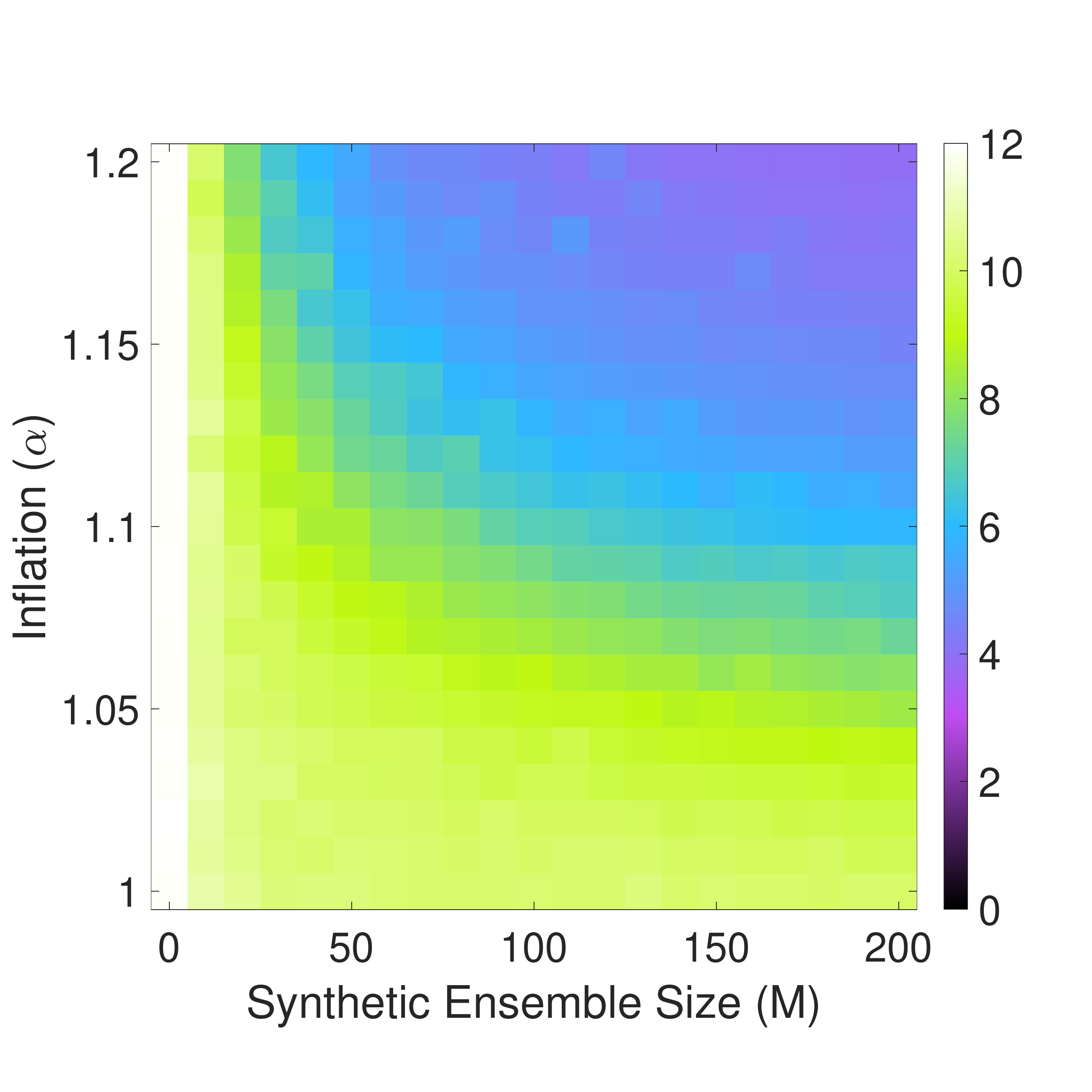}
        \caption{Analysis RMSE of the covariance shrinkage approach to particle rejuvenation (FETPF) for different values of the synthetic ensemble size $M$ and synthetic inflation factor $\alpha$, and for a dynamical ensemble size of $N=5$.}
        \label{fig:exp3}
    \end{figure}

    For all our experiments we use the Lorenz '63 system~\citep{lorenz1963deterministic}:
    \begin{equation}
        \begin{aligned}
            x' &= \sigma(y-x),\\
            y' &= x(\rho - z)-y,\\
            z' &= xy-\beta z,
        \end{aligned}		
    \end{equation}
    with canonical parameter values $\sigma = 10$, $\rho=28$, and $\beta=8/3$. We observe the first component, with Gaussian unbiased observation error, with a variance of $\*R=8$. The implementation is from our test problem suite~\citep{otpsoft,otp}.

    As discussed previously, the canonical choice for the shrinkage covariance is the identity matrix. It has been the authors' experience that for most dynamical systems the choice is poor. Moreover, the sequential data assimilation problem typically provides ways to calculate climatological approximations to the covariance. We take advantage of such techniques in this paper.
    
    The first type of climatological covariance that we investigate is that of the distribution over that of the whole manifold of the dynamics. The (trace-state normalized) matrix that is obtained by taking the temporal covariance of $50000$ sample points on the attractor of the canonical Lorenz '63 model is: 
    \begin{equation}\label{eqn:cov-P}
        \!P = \begin{bmatrix}[r]
            0.8616  &  0.8618  & -0.0148\\
            0.8618  &  1.1149  & -0.0035\\
            -0.0148  & -0.0035  &  1.0234
        \end{bmatrix},
    \end{equation}
    with condition number $15.88$.
    
    For testing multiple covariances, we run an ETPF with $N=100$ with 20000 evenly spaced samples over a time interval of $2400$ time units and calculate the trace-state normalized forecast covariances. Under a square Frobenius norm distance, we cluster the empirical covariance matrices of the same ensemble at different times using the $k$-means algorithm~\citep{karpatne2018book} into two clusters. The collection of climatological covariances for the Lorenz '63 thus consists of the centroids of each cluster, 
    \begin{equation}\label{eqn:cov-P1-P2}
        \begin{gathered}
            \!P_1 = \begin{bmatrix}[r]
                0.5017 &   0.5524 &  -0.4587\\
                0.5524  &  1.0731  & -0.6723\\
                -0.4587 &  -0.6723 &   1.4252
            \end{bmatrix}, \quad
            \!P_2 = \begin{bmatrix}
                0.5443  &  0.6830  &  0.4330\\
                0.6830  &  1.2748 &   0.6318\\
                0.4330   & 0.6318  &  1.1808
            \end{bmatrix},
        \end{gathered}
    \end{equation}
    with condition numbers $13.68$ and $16.98$ respectively.
    As can be seen, the clusters are mainly split by the correlation factors of $z$ with respect to the other variables being positive or negative. 
    
    In order to stay in line with other particle rejuvenation techniques, a heuristic that we use is inflation on the synthetic ensemble, so that it is possible to overcome deficiencies in its descriptive power:
    \begin{equation}\label{eqn:synthetic-ensemble-inflation}
        \!A^\|f \xleftarrow{} \alpha\!A^\|f,
    \end{equation}
    which is equivalent to assuming an inflated scaling factor $\mu$ in \eqref{eqn:mu}.
    We therefore have two parameters that we can configure in our rejuvenation technique: $M$, the size of our surrogate synthetic ensemble, and $\alpha$, the inflation applied to its realizations. It is the authors' experience that inflation should only be applied to the synthetic ensemble, and not the dynamical ensemble, otherwise the shrinkage estimate \eqref{eq:shr-rblw} does not lead to a stable algorithm.
    
    In our experiments we report the error of the analysis with respect to the truth, measured by the spatio-temporal RMSE, given by the formula,
    \begin{equation}
        \mathrm{RMSE}(\*x^{\|t},\, \bar{\*x}^{\|a}) = \sqrt{\frac{1}{n K} \sum_{i = 1}^T\sum_{j=1}^n \left(\left[\*x^{\|t}_i\right]_j - \left[\bar{\*x}^{\|a}_i\right]_j \right)^2 },
    \end{equation}
    where $T$ stands for the relevant measured timeframe of the experiments.
    
    Our first round of experiments reported in \Cref{fig:exp1} compares the canonical method of rejuvenation in  ETPF (and second order extensions) with the optimal rejuvenation factor of $\tau = 0.04$  in \eqref{eqn:canonicalrejuv} (computed by parameter search) to the stochastic covariance shrinkage technique for both Gaussian and Laplace samples. A dynamic ensemble size of $M=100$, the inflation factors $\alpha=\{1.0,1.2\}$. The target covariance \eqref{eqn:cov-P} is used. We perform  10000 assimilation steps, but discard the first $1000$ that are used for spinup. The time interval between successive observations  is $\Delta t = 0.12$. We perform $20$ independent runs and take the mean of the results to obtain an accurate estimate of the expected error. 
    
    Results in \Cref{fig:exp1} show that the stochastic covariance shrinkage technique converges to the same `correct' RMSE in the case of a large dynamic ensemble of $N=100$. Moreover, in the small ensemble case of $N=5$, methods that inflate the synthetic ensemble with factor $\alpha=1.2$ perform significantly better than those that do not. Laplace distributed synthetic samples also seem to slightly reduce the error.
    
    Our second round of experiments reported in \Cref{fig:exp2} makes use of multiple values of the climatological covariance $\!P$, namely those in \eqref{eqn:cov-P1-P2}. The rest of the setup is identical to the previous experiment. Again, for a large value of dynamic ensemble size $N=100$, the stochastic covariance shrinkage approach attains the correct error statistics. For a low ensemble size, however, there is virtually no difference between the Gaussian, Laplace, inflation, and no inflation stochastic covariance shrinkage methods as compared to the ETPF. 
    
    Our final round of experiments seeks to understand the effect of selecting the two free parameters, i.e., the synthetic ensemble size $M$ and the synthetic ensemble inflation factor $\alpha$. \Cref{fig:exp3} shows the spatio-temporal RMSE of a small dynamic ensemble ($N=5$) for various values of $M$ and $\alpha$, with Gaussian synthetic samples using the single target matrix \eqref{eqn:cov-P}. The results clearly show that in many operationally useful cases, it is necessary to have a sufficiently expressive synthetic ensemble, whose anomalies are sufficiently inflated.
    
    The results of the first two experiments show that adding additional synthetic information during assimilation is more effective than randomly perturbing the ensemble post-assimilation. The authors hypothesize that the results point strongly towards the need of intelligently, and adaptively choosing the target covariance matrices, and to the need for better operational calculation of the covariance shrinkage factor $\gamma$.

    %%%%%%%%%%%%%%%%%%%%%
    \section{Conclusions}
    \label{sec:conclusions}
    %%%%%%%%%%%%%%%%%%%%%

        This paper introduces a stochastic covariance shrinkage-based particle rejuvenation technique for the ensemble transport particle filter. Instead of reweighing existing prior information, the approach incorporates additional prior information into the ensemble through the use of synthetic anomalies. These anomalies come from climatological covariance information.
        Numerical experiments show that the use of climatological prior information to perform rejuvenation leads to reduced analyses errors for significantly smaller dynamical ensemble sizes than the original rejuvenation approach.
        
        We believe that the stochastic covariance shrinkage approach to importance sampling can be used not just for particle rejuvenation in the ETPF, but in other particle filters as well.

    %%%%%%%%%%%%%%%%%%%%%
    \section*{Acknowledgements}
    %%%%%%%%%%%%%%%%%%%%%
    This work was supported by awards NSF CDS\&E-MSS--1953113, DOE ASCR DE--SC0021313,
    and by the Computational Science Laboratory at Virginia Tech.

    %\section*{References}
    \bibliographystyle{copernicus}
    \bibliography{biblio}

\begin{thebibliography}{37}
\providecommand{\natexlab}[1]{#1}
\providecommand{\url}[1]{{\tt #1}}
\providecommand{\urlprefix}{URL }
\expandafter\ifx\csname urlstyle\endcsname\relax
  \providecommand{\doi}[1]{https://doi.org/\discretionary{}{}{}#1}\else
  \providecommand{\doi}{https://doi.org/\discretionary{}{}{}\begingroup
  \urlstyle{rm}\Url}\fi

\bibitem[{Acevedo et~al.(2017)Acevedo, de~Wiljes, and
  Reich}]{acevedo2017second}
Acevedo, W., de~Wiljes, J., and Reich, S.: Second-order accurate ensemble
  transform particle filters, SIAM Journal on Scientific Computing, 39,
  A1834--A1850, 2017.

\bibitem[{Anderson(2001)}]{anderson2001ensemble}
Anderson, J.~L.: An ensemble adjustment {K}alman filter for data assimilation,
  Monthly weather review, 129, 2884--2903, 2001.

\bibitem[{Anderson(2012)}]{anderson2012localization}
Anderson, J.~L.: Localization and sampling error correction in ensemble
  {K}alman filter data assimilation, Monthly Weather Review, 140, 2359--2371,
  2012.

\bibitem[{Asch et~al.(2016)Asch, Bocquet, and Nodet}]{asch2016data}
Asch, M., Bocquet, M., and Nodet, M.: {Data assimilation: methods, algorithms,
  and applications}, SIAM, 2016.

\bibitem[{Burgers et~al.(1998)Burgers, van Leeuwen, and
  Evensen}]{burgers1998analysis}
Burgers, G., van Leeuwen, P.~J., and Evensen, G.: Analysis scheme in the
  ensemble {K}alman filter, Monthly weather review, 126, 1719--1724, 1998.

\bibitem[{Chen et~al.(2009)Chen, Wiesel, and Hero}]{chen2009shrinkage}
Chen, Y., Wiesel, A., and Hero, A.~O.: Shrinkage estimation of high dimensional
  covariance matrices, in: 2009 IEEE International Conference on Acoustics,
  Speech and Signal Processing, pp. 2937--2940, IEEE, 2009.

\bibitem[{Chustagulprom et~al.(2016)Chustagulprom, Reich, and
  Reinhardt}]{chustagulprom2016hybrid}
Chustagulprom, N., Reich, S., and Reinhardt, M.: A hybrid ensemble transform
  particle filter for nonlinear and spatially extended dynamical systems,
  SIAM/ASA Journal on Uncertainty Quantification, 4, 592--608, 2016.

\bibitem[{{Computational Science Laboratory}(2020)}]{otpsoft}
{Computational Science Laboratory}: {ODE} Test Problems,
  \urlprefix\url{https://github.com/ComputationalScienceLaboratory/ODE-Test-Problems},
  2020.

\bibitem[{Evensen(1994)}]{evensen1994sequential}
Evensen, G.: {Sequential data assimilation with a nonlinear quasi-geostrophic
  model using Monte Carlo methods to forecast error statistics}, Journal of
  Geophysical Research: Oceans, 99, 10\,143--10\,162, 1994.

\bibitem[{Evensen(2009)}]{evensen2009data}
Evensen, G.: {Data assimilation: the ensemble Kalman filter}, Springer Science
  \& Business Media, 2009.

\bibitem[{Hunt et~al.(2007)Hunt, Kostelich, and Szunyogh}]{hunt2007efficient}
Hunt, B.~R., Kostelich, E.~J., and Szunyogh, I.: Efficient data assimilation
  for spatiotemporal chaos: A local ensemble transform {K}alman filter, Physica
  D: Nonlinear Phenomena, 230, 112--126, 2007.

\bibitem[{Jaynes(2003)}]{jaynes2003probability}
Jaynes, E.~T.: Probability theory: The logic of science, Cambridge university
  press, 2003.

\bibitem[{Kozubowski et~al.(2013)Kozubowski, Podg{\'o}rski, and
  Rychlik}]{kozubowski2013multivariate}
Kozubowski, T.~J., Podg{\'o}rski, K., and Rychlik, I.: Multivariate generalized
  {L}aplace distribution and related random fields, Journal of Multivariate
  Analysis, 113, 59--72, 2013.

\bibitem[{Lancewicki and Aladjem(2014)}]{lancewicki2014multi}
Lancewicki, T. and Aladjem, M.: Multi-target shrinkage estimation for
  covariance matrices, IEEE Transactions on Signal Processing, 62, 6380--6390,
  2014.

\bibitem[{Law et~al.(2015)Law, Stuart, and Zygalakis}]{law2015data}
Law, K., Stuart, A., and Zygalakis, K.: {Data assimilation: a mathematical
  introduction}, vol.~62, Springer, 2015.

\bibitem[{Liu(2008)}]{liu2008monte}
Liu, J.~S.: Monte Carlo strategies in scientific computing, Springer Science \&
  Business Media, 2008.

\bibitem[{Lorenc et~al.(2015)Lorenc, Bowler, Clayton, Pring, and
  Fairbairn}]{lorenc2015comparison}
Lorenc, A.~C., Bowler, N.~E., Clayton, A.~M., Pring, S.~R., and Fairbairn, D.:
  Comparison of hybrid-4DEnVar and hybrid-4DVar data assimilation methods for
  global NWP, Monthly Weather Review, 143, 212--229, 2015.

\bibitem[{Lorenz(1963)}]{lorenz1963deterministic}
Lorenz, E.~N.: Deterministic nonperiodic flow, Journal of the atmospheric
  sciences, 20, 130--141, 1963.

\bibitem[{McCann and Guillen(2011)}]{mccann2011five}
McCann, R.~J. and Guillen, N.: Five lectures on optimal transportation:
  geometry, regularity and applications, Analysis and geometry of metric
  measure spaces: lecture notes of the s{\'e}minaire de Math{\'e}matiques
  Sup{\'e}rieure (SMS) Montr{\'e}al, pp. 145--180, 2011.

\bibitem[{McCann et~al.(1995)}]{mccann1995existence}
McCann, R.~J. et~al.: Existence and uniqueness of monotone measure-preserving
  maps, Duke Mathematical Journal, 80, 309--324, 1995.

\bibitem[{Nino-Ruiz and Sandu(2015)}]{Sandu_2015_covarianceShrinkage}
Nino-Ruiz, E.~D. and Sandu, A.: Ensemble {Kalman} filter implementations based
  on shrinkage covariance matrix estimation, Ocean Dynamics, 65, 1423--1439,
  \doi{10.1007/s10236-015-0888-9}, 2015.

\bibitem[{Nino-Ruiz and Sandu(2017)}]{Sandu_2017_Covariance-Cholesky}
Nino-Ruiz, E.~D. and Sandu, A.: An ensemble {K}alman filter implementation
  based on modified {Cholesky} decomposition for inverse covariance matrix
  estimation, SIAM Journal on Scientific Computing, submitted,
  \urlprefix\url{https://arxiv.org/abs/1605.08875}, 2017.

\bibitem[{Nino-Ruiz and Sandu(2018)}]{Sandu_2019_Covariance-parallel}
Nino-Ruiz, E.~D. and Sandu, A.: Efficient parallel implementation of {DDDAS}
  inference using an ensemble {Kalman} filter with shrinkage covariance matrix
  estimation, Cluster Computing, pp. 1--11,
  \urlprefix\url{https://link.springer.com/article/10.1007/s10586-017-1407-1},
  2018.

\bibitem[{Nino-Ruiz et~al.(2015)Nino-Ruiz, Sandu, and Deng}]{Sandu_2015_SCALA}
Nino-Ruiz, E.~D., Sandu, A., and Deng, X.: A parallel ensemble {Kalman} filter
  implementation based on modified {Cholesky} decomposition, in: Proceedings of
  the 6th Workshop on Latest Advances in Scalable Algorithms for Large-Scale
  Systems, vol. Supercomputing 2015 of {\em ScalA '15\/}, Austin, Texas,
  \doi{10.1145/2832080.2832084}, 2015.

\bibitem[{Olver et~al.(2010)Olver, Lozier, Boisvert, and Clark}]{olver2010nist}
Olver, F.~W., Lozier, D.~W., Boisvert, R.~F., and Clark, C.~W.: NIST handbook
  of mathematical functions hardback and CD-ROM, Cambridge university press,
  2010.

\bibitem[{Petrie(2008)}]{petrie2008localization}
Petrie, R.: Localization in the ensemble {K}alman filter, MSc Atmosphere, Ocean
  and Climate University of Reading, 2008.

\bibitem[{Popov and Sandu(2020)}]{popov2020explicit}
Popov, A.~A. and Sandu, A.: An Explicit Probabilistic Derivation of Inflation
  in a Scalar Ensemble {K}alman Filter for Finite Step, Finite Ensemble
  Convergence, arXiv preprint arXiv:2003.13162, 2020.

\bibitem[{Popov et~al.(2020)Popov, Sandu, Nino-Ruiz, and
  Evensen}]{popov2020stochastic}
Popov, A.~A., Sandu, A., Nino-Ruiz, E.~D., and Evensen, G.: A Stochastic
  Covariance Shrinkage Approach in Ensemble Transform {K}alman Filtering, 2020.

\bibitem[{Rao et~al.(2017)Rao, Sandu, Ng, and Nino-Ruiz}]{rao2017robust}
Rao, V., Sandu, A., Ng, M., and Nino-Ruiz, E.~D.: {Robust Data Assimilation
  using $L_1$ and Huber norms}, SIAM Journal on Scientific Computing, 39,
  B548--B570, 2017.

\bibitem[{Reich(2013)}]{reich2013nonparametric}
Reich, S.: A nonparametric ensemble transform method for Bayesian inference,
  SIAM Journal on Scientific Computing, 35, A2013--A2024, 2013.

\bibitem[{Reich and Cotter(2015)}]{reich2015probabilistic}
Reich, S. and Cotter, C.: {Probabilistic forecasting and Bayesian data
  assimilation}, Cambridge University Press, 2015.

\bibitem[{Roberts et~al.(2019)Roberts, Popov, and Sandu}]{otp}
Roberts, S., Popov, A.~A., and Sandu, A.: {ODE} Test Problems: a {MATLAB} suite
  of initial value problems, CoRR, abs/1901.04098,
  \urlprefix\url{http://arxiv.org/abs/1901.04098}, 2019.

\bibitem[{Ruiz et~al.(2014)Ruiz, Sandu, and Anderson}]{Sandu_2014_EnKF_SMF}
Ruiz, E.~N., Sandu, A., and Anderson, J.: An efficient implementation of the
  ensemble {K}alman filter based on an iterative {Sh}erman-{M}orrison formula,
  Statistics and Computing, pp. 1--17, \doi{10.1007/s11222-014-9454-4}, 2014.

\bibitem[{Strogatz(2018)}]{strogatz2018nonlinear}
Strogatz, S.~H.: Nonlinear dynamics and chaos with student solutions manual:
  With applications to physics, biology, chemistry, and engineering, CRC press,
  2018.

\bibitem[{Tan et~al.(2018)Tan, Steinbach, Karpatne, and
  Kumar}]{karpatne2018book}
Tan, P.-N., Steinbach, M., Karpatne, A., and Kumar, V.: Introduction to Data
  Mining (2nd Edition), Pearson, 2nd edn., 2018.

\bibitem[{Villani(2003)}]{villani2003topics}
Villani, C.: Topics in optimal transportation, 58, American Mathematical Soc.,
  2003.

\bibitem[{Zhang et~al.(2010)Zhang, Liu, and Oliver}]{zhang2010ensemble}
Zhang, Y., Liu, N., and Oliver, D.~S.: Ensemble filter methods with perturbed
  observations applied to nonlinear problems, Computational Geosciences, 14,
  249--261, 2010.

\end{thebibliography}

\end{document}